\documentclass{svproc}

\usepackage{amsfonts}
\usepackage{amssymb}
\usepackage{amsmath}
\usepackage{stmaryrd}

 \def\cb{{\cal B}} 
  
\def\cg{{\cal G}} \def\ch{{\cal H}} 
 \def\ck{{\cal K}} 
  
\def\cp{{\cal P}}   \def\cs{{\cal S}}

\def\bbc{\mathbb{C}}
\def\bbr{\mathbb{R}}

\def\eqn#1{\begin{equation}\label{#1}}
\def\ee{\end{equation}}

\def\bea{\begin{eqnarray}}
\def\eea{\end{eqnarray}}

\def\eqnn#1{\begin{eqnarray}\label{#1}}
\newcommand{\eqna}[1]{\begin{subequations} \label{#1}
\begin{eqnarray}}
\def\eena{\end{eqnarray}
\end{subequations}}

\def\k{\kappa}

\def\nn{\nonumber}
\def\ha{{\textstyle{1\over2}}}
\def\d{\delta}
\def\L{\Lambda} \def\l{\lambda}

\def\rf#1#2{(\ref{#1}{#2})}

\def\trq{{\textstyle{3\over4}}}
\def\frq{{\textstyle{5\over4}}}

   \def\trh{{\textstyle{\frac{3}{2}}}}
   \def\nt{\noindent}
\def\qrt{{\textstyle{1\over4}}}

\usepackage{url}

\begin{document}
\mainmatter

\title{On Reducible Verma Modules over \\[3pt]
Jacobi Algebra}
\titlerunning{Jacobi Algebra}

\author{V.K. Dobrev} 
\institute{Institute of Nuclear Research and Nuclear Energy,\\
 Bulgarian Academy of Sciences, \\
72 Tsarigradsko Chaussee, 1784 Sofia, Bulgaria}

 \maketitle              

\begin{abstract}
With this paper we start the study of reducible representations of the Jacobi algebra
with the ultimate goal of constructing differential operators invariant w.r.t. the Jacobi algebra.
In this first paper we show examples of the low level singular vectors of Verma modules
over the Jacobi algebra. According to our methodology these will produce the invariant differential operators.
\keywords{Jacobi algebra, Verma modules, singular vectors}
\end{abstract}


\section{Introduction}

The role of nonrelativistic symmetries in  theoretical physics was
always important.  Currently one of the  most popular fields in
theoretical physics - string theory, pretending to be a universal
theory - encompasses together relativistic quantum field theory,
classical gravity, and certainly, nonrelativistic quantum mechanics,
in such a way that it is not even necessary to separate these
components.

Since the cornerstone of quantum mechanics is the Schr\"o\-din\-ger
equation then
 it is not a surprise that the Schr\"o\-din\-ger group -
  the group that is the maximal group of
symmetry of the Schr\"o\-din\-ger equation - was the first to play
a prominent role in theoretical physics.
The latter is natural since originally the Schr\"odinger group,
actually the Schr\"odinger algebra, was introduced  in
\cite{Nie,Hag} as a nonrelativistic limit of
the vector-field realization of the conformal algebra.
For a review on these developments we refer to \cite{VKD4}.

Another interesting non-relativistic example is the Jacobi algebra \cite{EiZa,BeSc} which is the
semi-direct sum of the Heisenberg algebra and the $sp(n)$ algebra. Actually the lowest case of the
Jacobi algebra coincides with the lowest case of the Schr\"odinger algebra which makes it interesting
to apply to the Jacobi algebra the methods we applied to the Schr\"odinger algebra. This is a project we start
in the present short paper. Actually here we give as examples the low level singular vectors of Verma modules
over the Jacobi algebra.

\section{Preliminaries}

The Jacobi algebra is the semi-direct sum $\cg_n:= \ch_n\niplus 
sp(n,\bbr)_{\bbc}$ \cite{EiZa,BeSc}. The Heisenberg algebra
$\ch_n$ is generated by the boson creation (respectively,
annihilation) operators~${a}_i^{+}$~(${a}^-_i$),~$i,j
=1,\dots,n$, which verify the canonical commutation relations
\eqn{heis} \big[a^-_i,a^{+}_j\big]=\delta_{ij}, \qquad [a^-_i,a^-_j]
= \big[a_i^{+},a_j^{+}\big]= 0 . \ee
$\ch_n$ is an
ideal in $\cg_n$, i.e., $[\ch_n,\cg_n]=\ch_n$,
determined by the commutation relations (following the notation of \cite{Berc}):
\eqna{haspn}
&&\big[a^{+}_k,K^+_{ij}\big] = [a^-_k,K^-_{ij}]=0, \\
&&{} [a^-_i,K^+_{kj}] = \tfrac{1}{2}\delta_{ik}a^{+}_j+\tfrac{1}{2}\delta_{ij}a^{+}_k ,\qquad
 \big[K^-_{kj},a^{+}_i\big] = \tfrac{1}{2}\delta_{ik}a^-_j+\tfrac{1}{2}\delta_{ij}a^-_k , \\
&& \big[K^0_{ij},a^{+}_k\big] = \tfrac{1}{2}\delta_{jk}a^{+}_i,\qquad
\big[a^-_k,K^0_{ij}\big]= \tfrac{1}{2}\delta_{ik}a^-_{j} .
\eena
 $K^{\pm,0}_{ij}$ are the generators of the $\cs_n ~\equiv~ sp(n,\bbr)_{\bbc}$ algebra:
\eqna{baspn}
&& [K_{ij}^-,K_{kl}^-] = [K_{ij}^+,K_{kl}^+]=0 , \qquad 2\big[K^-_{ij},K^0_{kl}\big] = K_{il}^-\delta_{kj}+K^-_{jl}\delta_{ki}\label{baza23}, \\
&& 2[K_{ij}^-,K_{kl}^+] = K^0_{kj}\delta_{li}+
K^0_{lj}\delta_{ki}+K^0_{ki}\delta_{lj}+K^0_{li}\delta_{kj}\\
&& 2\big[K^+_{ij},K^0_{kl}\big] = -K^+_{ik}\delta_{jl}-K^+_{jk}\delta_{li},\quad
 2\big[K^0_{ji},K^0_{kl}\big] = K^0_{jl}\delta_{ki}-K^0_{ki}\delta_{lj} . 
\eena

In order to implement our approach we introduce a triangular decomposition of ~$\cg_n$~:
\eqn{deca} \cg_n\,=\, \cg_n^+\oplus \ck_n\oplus \cg_n^- \ ,\ee
using the triangular decomposition ~$\cs_n ~=~ \cs_n^+\oplus \ck_n\oplus \cs_n^- $,
where:
\eqnn{decpm} && \cg_n^\pm ~=~ \ch_n^\pm \oplus \cs_n^\pm \\
&& \ch^\pm_n ~=~ {\rm l.s.} \{\, {a}_i^{\pm} : i=1,\dots,n\}\ , \nn\\
&& \cs_n^+ ~=~ {\rm l.s.} \{\, K^+_{ij} ~:~ 1\leq i\leq j\leq n \}
\oplus {\rm l.s.} \{\, K^0_{ij} ~:~ 1\leq i<  j\leq n \} \nn\\
&& \cs_n^- ~=~ {\rm l.s.} \{\, K^-_{ij} ~:~ 1\leq i\leq j\leq n \}
\oplus {\rm l.s.} \{\, K^0_{ij} ~:~ 1\leq j < i \leq n \} \nn\\
&&\ck_n ~=~ {\rm l.s.} \{\, K^0_{ii} ~:~ 1\leq  i \leq n \}
\nn\eea
Note that the subalgebra ~$\ck_n$~ is abelian and
 is a Cartan subalgebra of ~$\cs_n$. Furthermore, not only ~$\cs_n^\pm$, but also
 ~$\cg_n^\pm$~ are its eigenspaces:
\eqn{ckcg} [ \ck_n, \cg_n^\pm] ~=~ \cg_n^\pm \ee
Thus, ~$\ck_n$~ plays for ~$\cg_n$~ the role that
Cartan subalgebras are playing for semi-simple Lie algebras.

\section{Case $\cg_2$}

Note that the algebra ~$\cg_1$~ is isomorphic to the (1+1)-dimensional Schr\"odinger algebra (without
central extension). The representations of the latter are well known, cf. \cite{DDM,AiDo,DLMZ,VKD4}.
Thus, we study the first new  case of the ~$\cg_n$~ series, namely, ~$\cg_2$.

For simplicity, we introduce the following notations for the basis of ~$\cs_2$~:
\eqna{bas2} \cs^+ ~:~&&~ b^+_i ~\equiv~ K^+_{ii}\ , ~~i=1,2; \quad
c^+ ~\equiv~ K^+_{12}\ ,\quad d^+ ~\equiv K^0_{12} \\
\cs^- ~:~&&~ b^-_i ~\equiv~ K^-_{ii}\ , ~~i=1,2; \quad
c^- ~\equiv~ K^-_{12}\ ,\quad d^- ~\equiv K^0_{21} \\
\ck ~:~&&~ h_i ~\equiv~ K^0_{ii} \ , ~~i=1,2.
\eena
Next, using \eqref{haspn} and \eqref{baspn} we give the eigenvalues of the basis of
~$\cg^+$~ w.r.t. ~$\ck$~:
\eqnn{eig2}
&& h_1 ~:~ (b^+_1, b^+_2, c^+, d^+ , a_1^+, a_2^+) ~:~ (1,0,\ha,\ha,\ha,0) \ , \\
&& h_2 ~:~ (b^+_1, b^+_2, c^+, d^+ , a_1^+, a_2^+) ~:~ (0,1,\ha,-\ha,0,\ha) \ , \nn\eea
(e.g., ~$[h_1,b^+_1] = b^+_1$, ~~$[h_2,d^+] = -\ha d^+$, etc).
Naturally, the  eigenvalues of the basis of ~$\cg^-$~ w.r.t. ~$\ck$~ are obtained from \eqref{eig2}
by multiplying every eigenvalue by (-1).

Next we introduce the following grading of the basis of ~$\cg^+_2$:
\eqn{grad}
(b^+_1, b^+_2, c^+, d^+ , a_1^+, a_2^+) ~:~ (2\d_1,2\d_2,\d_1+\d_2,\d_1-\d_2,\d_1,\d_2)
\ee
The grading of the   ~$\cs^+_2$~ part of the basis follows from the root system of  ~$\cs^+_2$,
while the grading of the ~$\ch^+_2$~ part of the basis is determined by consistency with commutation
relations \eqref{haspn}. It is consistent also with  formulae \eqref{eig2}.

Naturally, the grading  of the basis of ~$\cg^-$~ w.r.t.  are obtained from \eqref{grad}
by multiplying every grading  by (-1).

\section{Verma modules and singular vectors}

\subsection{Definitions}

We shall introduce Verma modules over the Jacobi algebra analogously to the case of
  of semi-simple algebras.
Thus, we define a  lowest weight ~{\it Verma module} ~$V^\L$~ over $\cg_n$   as
the lowest  weight module over ~$\cg_n$~ with lowest  weight ~$\L \in \ck_n^*$~ and
lowest weight vector ~$v_0 \in V^\L$, induced from the
one-dimensional representation ~$V_0 \cong \bbc v_0$~ of
~$U(\cb_n)$~, (where ~$\cb_n  = \ck_n \oplus \cg_n^-$~ is a Borel
subalgebra of ~$\cg_n$), such that:
\eqnn{indb}
& &X ~v_0 ~~=~~ 0 , \quad \forall\, X\in \cg_n^- \cr
&&H ~v_0 ~~=~~ \L(H)~v_0\,, \quad \forall\, H \in \ck_n \eea

Pursuing the analogy with the semi-simple case and following our approach we are interested
in the cases when the Verma modules are reducible. Namely, we are interested
in the cases when a Verma module $V^\L$  contains an invariant submodule
which is also a Verma module ~$V^{\L'}$, where ~$\L'\neq \L$, and holds the analog of
\eqna{indbb}
& &X ~v'_0 ~~=~~ 0 , \quad \forall\, X\in \cg_n^- \\
&&H ~v'_0 ~~=~~ \L'(H)~v'_0\,, \quad \forall\, H \in \ck_n \eena
Since ~$V^{\L'}$~ is an invariant submodule then there should be a mapping such that
~$v'_0$~ is mapped to a singular vector ~$v_s \in V^\L$~ fulfilling exactly \eqref{indbb}.
Thus, as in the semi-simple case there should be a polynomial ~$\cp$~ of ~$\cg_n^+$~ elements
which is eigenvector of $\ck_n$: ~$[H,\cp] = \L'(H)\cp$, ($\forall H\in\ck_n$), and then we would have: ~$v_s ~=~\cp v_0\,$.

\subsection{Case $\cg_2$}

We shall consider several examples of reducible Verma modules with different weights.

\subsubsection{Weight $2\d_1$}

As first example we try to find a singular vector of weight ~$\L' \sim 2\d_1\,$.  There are six possible terms
in ~$U(\cg_2)$~ with this weight, thus, we try:
\eqn{delta11}
v_s^{2\d_1} ~=~ \big( \nu_1 b^+_1 + \nu_2 c^+d^+ + \nu_3 b^+_2 (d^+)^2 + \nu_4 (a_1^+)^2
+ \nu_5 a_1^+ a_2^+ d^+ + \nu_6 (a_2^+)^2 (d^+)^2 \big) v_0
\ee
where ~$\nu_k$~ are numerical coefficients which may be fixed when we impose
\rf{indbb}{a} on \eqref{delta11}. (Note that \rf{indbb}{b} is fulfilled by every term of \eqref{delta11}.)

After we impose \rf{indbb}{a} on \eqref{delta11} we find the solution:
\eqnn{sol2d} &&\L(H_1) ~=~ \trq\,, \quad ~\nu_3 ~=~ -2\nu_6,\nn\\
&& \nu_1~=~ - \nu_6(\L(H_2) - \L(H_1))( 2\L(H_2) - 2\L(H_1) -1),\nn\\
&& \nu_2 ~=~  2\nu_6 ( 2\L(H_2) - 2\L(H_1) -1), \nn\\
&& \nu_4 ~=~  \nu_6(\L(H_2) - \L(H_1))( \L(H_2) - \L(H_1) -\ha),\nn\\
&& \nu_5 ~=~  -\nu_6 ( 2\L(H_2) - 2\L(H_1) -1).\eea
Thus, the singular vector is:
\eqnn{singb} v_s^{2\d_1} ~&=&~ \nu_6 \big(\
  (\L(H_2) - \trq )( \L(H_2) - \frq  ) ((a_1^+)^2 -2 b^+_1)\ +~ \nn\\
&&+~ 2 ( \L(H_2) - \frq) (2c^+ -   a_1^+ a_2^+) d^+ \ + \nn\\
&&+~  ((a_2^+)^2  -2 b^+_2 )(d^+)^2\ \big)v_0\  , \quad\
\L(H_1) ~=~ \trq
\eea

\subsubsection{Weight $2\d_2$}

As next example we try to find a singular vector of weight ~$\L' \sim 2\d_2\,$.  The
possible singular vector is:
\eqn{delta12}
v_s^{2\d_2} ~=~ \big(\ \mu_1 b^+_2 +   \mu_2 (a_2^+)^2\ \big) v_0
\ee
Imposing \rf{indbb}{a} on \eqref{delta12} we obtain:
\eqn{so2d} \L(H_2) ~=~ \qrt\,, \quad ~\mu_1 ~=~ -2\mu_2, \ee
Thus, the singular vector is:
\eqn{sing2} v_s^{2\d_2} ~=~ \mu_2 ((a_2^+)^2 - 2 b^+_2 )v_0 \ , \quad
\L(H_2) ~=~ \qrt
 \ee

\subsubsection{Weight $\d_1+\d_2$}

Next we try a singular vector of weight ~$\L' \sim \d_1+\d_2\,$.  The
possible singular vector is:
\eqn{delt12}
v_s^{\d_1+\d_2} ~=~ \big(\  \k_1 c^+ + \k_2 b^+_2 d^+
+ \k_3 a_1^+ a_2^+  + \k_4 (a_2^+)^2 d^+\ \big) v_0
\ee
Imposing \rf{indbb}{a} on \eqref{delt12} we obtain:
\eqnn{sod12} && \L(H_2) ~=~ \trh - \L(H_1)\,, \quad \k_1 = (3-4h(1))\k_4\ , \nn\\
&& \k_2 = -2\k_4\ , \quad \k_3 = (2h(1) - \trh)\k_4 \eea
Thus, the singular vector is:
\eqn{sing3} v_s^{\d_1+\d_2} ~=~ \k_4 \big(\ (\trh-2h(1)) (2c^+ - a_1^+ a_2^+)  +  ((a_2^+)^2 -2b^+_2) d^+ \ \big) v_0 \ee

\subsubsection{Weight $\d_1-\d_2$}

Next we try a singular vector of weight ~$\L' \sim \d_1-\d_2\,$.  The
only possible singular vector is:
\eqn{delm12}
v_s^{\d_1-\d_2} ~=~ \l  d^+ v_0
\ee
Imposing \rf{indbb}{a} on \eqref{delm12} we obtain that ~$v_s^{\d_1-\d_2}$~ is a singular vector iff:
\eqn{som12}  \L(H_2) ~=~  \L(H_1)\ee

\subsubsection{Weight $\d_1$}

Next we try a singular vector of weight ~$\L' \sim \d_1\,$. The
possible singular vector is:
\eqn{deld1}
v_s^{\d_1} ~=~  \big(\ \l_1 a^+_1 +   \l_2 a_2^+ d^+\ \big) v_0
\ee
Imposing \rf{indbb}{a} on \eqref{deld1} we obtain:
\eqn{resd1} \l_1 = \l_2 =0 \ee
Thus, there is no singular vector of weight $\d_1\,$.

\subsubsection{Weight $\d_2$}

Finally,  we try a singular vector of weight ~$\L' \sim \d_2\,$.  The
only possible singular vector is:
\eqn{deld2}
v_s^{\d_2} ~=~ \mu a_2^+ v_0
\ee
Imposing \rf{indbb}{a} on \eqref{deld2} we obtain:
\eqn{resd1} \mu =0 \ee
Thus, there is no singular vector of weight $\d_2\,$.

\subsubsection{Weight $3\d_2$}

 The only possible singular vector is:
\eqn{del3d2}
v_s^{\d_2} ~=~ \mu b_2^+ a_2^+ v_0 + \nu (a_2^+)^3 v_0
\ee
Imposing \rf{indbb}{a} on \eqref{del3d2} we obtain:
\eqn{res3d2}  \mu =\nu= 0  \ee
Thus, there is no singular vector of weight $3\d_2\,$.

\bigskip

\section*{Acknowledgments}

\nt
The author acknowledges partial support from Bulgarian NSF Grant DN-18/1.

\end{document}